\documentclass[a4paper,11pt]{article}

\usepackage{geometry,enumerate,amsmath,amssymb,latexsym,theorem,QED,epic,eepic}

\usepackage[all]{xy}
\usepackage[ps2pdf, bookmarks=true,
            bookmarksnumbered=true, breaklinks=true,
            pdfstartview=FitH, hyperfigures=false,
            plainpages=false, naturalnames=true,
            colorlinks=true,pagebackref=true,
            pdfpagelabels]{hyperref}
\usepackage[usenames]{color}

{\theorembodyfont{\rmfamily}\newtheorem{example}{Example}[section]}
{\theorembodyfont{\rmfamily}\newtheorem{defn}[example]{Definition}}

\newtheorem{thm}[example]{Theorem}
{\theorembodyfont{\rmfamily}\newtheorem{rem}[example]{Remark}}

\newcommand{\xdirects}[2]{\def\objectstyle{\scriptstyle} \objectmargin={0pt}
\xy
(0,0)*+{}="a",(0,-6)*+{\rule{0em}{1.5ex}#2}="b",(7,0)*+{\;#1}="c"
\ar@{->} "a";"b" \ar @{->}"a";"c" \endxy }

\newcommand{\T}{\textstyle}
\def\s{\mathsf s}
\def\t{\mathsf t}
\def\Ob{\operatorname{Ob}}
\newcommand{\brclass}[1]{\langle\!\langle #1 \rangle \! \rangle}
\def\xybiglabels{\def\labelstyle{\textstyle}}

\newcommand{\bl}{\mbox{\rule{0.08em}{1.7ex}\hspace{-0.00em}\rule{0.7em}{0.2ex}}}

\newcommand{\br}{\mbox{\rule{0.7em}{0.2ex}\hspace{-0.04em}\rule{0.08em}{1.7ex}}}

\newcommand{\tr}{\mbox{\rule[1.5ex]{0.7em}{0.2ex}\hspace{-0.03em}\rule{0.08em}{1.7ex}}}

\newcommand{\tl}{\mbox{\rule{0.08em}{1.7ex}\rule[1.54ex]{0.7em}{0.2ex}}}

\newcommand{\hh}{\mbox{\rule{0.7em}{0.2ex}\hspace{-0.7em}\rule[1.5ex]{0.70em}{0.2ex}}}

\newcommand{\vv}{\mbox{\rule{0.08em}{1.7ex}\hspace{0.6em}\rule{0.08em}{1.7ex}}}

\begin{document}
\title{`Double modules', double categories and groupoids,\\ and a new homotopical double groupoid}
\author{Ronald Brown\thanks{email: r.brown@bangor.ac.uk. This research
was partially supported by INTAS grant 93-436 ext `Algebraic
K-theory, groups and categories', and a Leverhulme Emeritus
Fellowship (2002-2004). }
\\ School of Computer Science
\\   Bangor University\\ Gwynedd LL57 1UT, U.K. }
\maketitle

\begin{center}
  Bangor Math Preprint 0901
\end{center}

\begin{abstract}
We give a rather general construction of double categories and so
double groupoids from a structure we call a `double module'.

We also give a homotopical construction of a double groupoid from a
triad consisting of a space, two subspaces, and a set of base
points, under a condition which also implies that this double
groupoid contains two second relative homotopy groups.\footnote{AMSCLass2000: 18B40,18D05,55E30 \\
KEYWORDS: double categories, double groupoids, crossed modules}
\end{abstract}
\section*{Introduction}

Double categories were introduced by Ehresmann \cite{Ehresmann-65}
as an example of structured category. Examples of double groupoids
were shown to arise from crossed modules in \cite{BS1,BS76} and this
result was applied in \cite{BH78:sec} to give a 2-dimensional
version of the Van Kampen Theorem for the fundamental group, namely
a colimit theorem for the fundamental crossed module of a based
pair. These results were generalised to crossed modules of
groupoids, and  to higher dimensions,  in \cite{BH81:col}. A more
general version of the construction of double groupoids was given in
\cite{Br-Ma} in terms of a {\em core diagram}. Another construction
of double groupoids due to the author was taken up by Lu and
Weinstein in \cite{Lu-weinstein} for purposes of Poisson groupoids.

The first aim of this note is to give a generalisation of this last
construction. We do not obtain an equivalence of categories, and so
in effect the construction shows that double groupoids can be quite
complex objects. Perhaps they should be considered among the basic
structures in mathematics. A classification of double groupoids is
given in \cite{and-natale-structure}.

The second aim (Section \ref{sec:newhomotopdg}) is to give a new
construction of a double groupoid for a topological space $X$ with
three subspaces $A,B,C$ such that $C \subseteq A \cap B$. Here $C$
is thought of as a set of base points, and the double groupoid is
well defined if the two induced morphisms $\pi_2(A,c) \to
\pi_2(X,c), \; \pi_2(B,c) \to \pi_2(X,c)$ have the same image. Under
this condition, the double groupoid also contains the two relative
homotopy groups $\pi_2(X,A,c), \; \pi_2(X,B,c)$ for all $c \in C$.
This  is an extension of results of \cite{BH78:sec}. Thus this
construction has the advantage of generality, symmetry, and multiple
compositions in either directions, advantages not available for the
traditional relative homotopy groups.

The relation between the two constructions is unclear, and it is
hoped that this paper will encourage further study of the area.

The bibliography gives other uses of double categories and
groupoids,
\cite{S,BrownMosaTAC,dawson-pare-assoc,ehresmann-catdoublestruct,BJ04,mack-sympldouble},
but is not intended to be exhaustive.

\section{Double modules and double categories}

We are trying to find a mathematical expression for the following
diagram and associated equations:
\begin{equation} \label{4square} \def\labelstyle{\textstyle}
\vcenter{\xymatrix@M=0pt @=3pc{ \ar @{-} [rr] ^(0.25)a ^(0.75)b
\ar @{-}[dd] _(0.25)u _(0.75)x  \ar@{} [dr] |m & \ar @{-}[dd]
|(0.25)v |(0.75)
y  \ar @{} [dr]| n & \ar @{-}[dd] ^(0.25)w ^(0.75)z \\
\ar@{-} [rr] |(0.25) c |(0.75) d \ar @{} [dr]| p & \ar@{} [dr]| q& \\
\ar@{-} [rr] _(0.25)e _(0.75) f & & }}  \qquad \xdirects{2}{1}
\end{equation}

We interpret each square as a 2-cell: thus $m$ is thought of as
\begin{align*}
av &\Leftarrow uc: m\\
\intertext{This might be formalised as the equation} av &=ucm.
\end{align*}   So our four squares give us {\em boundary equations}
\begin{alignat*}{2} av&=ucm, \qquad &
bw &= vdn, \\
cy &= xep, & dz&= yfq.
\end{alignat*}
We would like to `compose' such squares. The base point of each
square is thought of as the bottom right hand corner, so that is
where the centre values $m,n,p,q$ are `located'. So in order to
`compose' $m$ and $n$ we need to translate $m$ to the same corner as
$n$. Thus we assume an action $(m,d) \mapsto m^d$ satisfying $md=
dm^d$, and similarly $my=ym^y$. (The data and axioms for this will
be explained below: here we are concerned with formulae!)

Thus we deduce from the above rules that
\begin{alignat*}{2}
abw &=avdn  \qquad &avy&= ucm y \\
&= ucmdn  \qquad &&= ucym^y  \\
&= ucd m^d n \qquad&  &= uxep m^y.\\
\end{alignat*}

 We now construct from the above `data' a double category  $D$.
 The  squares  of $D$  will be quintuples
$ (m:u {\;}^{\T a} _{\T c}\;  v)$ such that  $ av=ucm  $ (compare
\cite{ehresmann-quintettes,S,BrownMosaTAC}).

The horizontal and vertical compositions are given by:
\begin{alignat*}{2}  (m:u {\;} ^{\textstyle a} _{\T c} {\;} v) &\circ _2
(n:v {\;}^{\T b} _{\T d} \; w) &&= (m^d \,n:u {\;\, } ^{\T ab}
_{\T cd } \;\, w) ,
\\
   (m:u {\;} ^{\T a} _{\T c}
\;  v) &\circ _1 (p:x {\;} ^{\T c}_{\T e}\; y) &&= (p\,m^y:ux{\;}
^{\T a}_{\T e}\; vy) . \end{alignat*}

In virtue of the above calculations, given the data of diagram
\eqref{4square}, these give compositions with the correct boundary
equations.

We also would like the interchange law, namely that the two possible
ways of composing the four squares in diagram \eqref{4square} give
the same answer. A direct calculation shows that this is equivalent
to the rule
\begin{equation}
  m^{yf} q = q m^{dz},
\end{equation}
again given the boundary equation $dz= yfq$, and this  makes
geometric sense in terms of the choices of `transporting' $m$ to the
bottom right hand corner of the square involving $q$.

\vspace{2ex}

Now we need to give a  structure in which the above calculations
make sense.

Our notation for categories will be that we write $\s,\t: C \to
\Ob{C}$ for the source and target maps of the category $C$ so that
an arrow $c$ of $C$ is an arrow $c: \s c \to \t c$, and
composition $cd$ is defined if and only if $\t c = \s d$.

\begin{defn} { A   {\em double $P$-module} consists  of  three  morphisms  of
categories all over the identity on objects: $$ \xymatrix{  M \ar
[dr] ^{\mu} & H \ar [d] ^{\phi} \\ V \ar [r] _{\psi} & P } $$ such
that $M$  is totally  intransitive, i.e. is a union of monoids, so
that $\s =\t$ on $M$. We write $M(x)$ for $M(x,x)$. Further, there
are given right actions of both $H$ and $V$ on $M$. This means that
if $m \in M(x)$ and $d \in H(x,y)$ then there is defined $m^d \in
M(y)$ and the usual axioms for an action are satisfied, namely
$(mm_1)^h=m^hm_1^h $, $m^{hk}=(m^h)^k$,  whenever these make sense,
and $1^h=1, m^1=m$, and similarly for the action of $V$ on $M$.

We do not suppose these actions commute, but nonetheless we agree
to write $m^{dz}$ when $\t m = \s d$ and $\t d = \s z$. However
$dz$ is here interpreted formally, or, if you like, as an arrow
 of the free product category $H * V$ over the same set of objects as
 $H,V,M$. Thus $\t m^{dz}= \t z$.

In order to write our axioms in a way which agrees with the above
calculations, we agree  `evaluate in $P$' means apply the
morphisms $\mu, \phi, \psi$ to the given equation to give an
equation in $P$. Thus the equation `$ucm =av$ evaluated in $P$'
means not only that
$$ (\psi u)(\phi c) (\mu m) = (\psi a)( \phi v)$$
but also that the equation makes sense in that $$\s a=\s u, \t a =
\s v, \t u= \s c , \t c =\t v \t m. $$ Similarly, $$md=dm^d \mbox{
evaluated in $P$ }$$  means that $ \t m = \s d$ and$$(\mu m )
(\phi d) = (\phi d) (\mu m^d).
$$

With this agreed, the axioms are:  if $m ,q\in M, d,f \in H, y,z
\in V$
\begin{enumerate}[(i)] \item then $$md= dm^d, my = y m^y
, \mbox{ evaluated in $P$ };$$  \item if also $yfq=dz$ in $P$, then
in $M$,
$$ m^{yf}q= q m^{dz}.
$$
\end{enumerate}}\qed
\end{defn}

Now we have our main result, which follows from what was written
above:
\begin{thm}
Given a double $P$-module as above, then the compositions
$\circ_1, \circ _2$ give the structure of double category, which
is a double groupoid if all of $M,H,V,P$ are groupoids.
\end{thm}

     Some special cases are of interest.
\begin{example}{   If  $M$  consists only of identities,  then   $D$
is the  double groupoid of squares from  $H$  and  $V$  which
``commute in $P$". This is used in \cite{Lu-weinstein}.} \qed
\end{example}
\begin{example}{   Suppose that $ V = P $ and  $\psi$    is
 the identity. Then
we obtain a diagram $$\xymatrix { M \ar [r] ^{\mu} &  P & \ar [l]
_{\phi} H} $$
 together with actions of  $P$
and  $H$  on  $M$.  The axioms now imply that  $\mu$ is a crossed
module.} \qed
\end{example}
\begin{example}{  Let  $H$  and  $V$  be subgroups of the group  $P$ and
 let $M$ be a
subgroup of  $P$  which  is  normal  in  both $H$   and $V$. Then
the inclusions, and the conjugation actions of  $H$  and $V$ on $M$
give the structure of a double   $P$-module, from which we can
obtain a double groupoid.} \qed
\end{example}
\begin{example}{  A semicore  diagram  is defined in \cite{Br-Ma}  to consist of a
commutative diagram of morphisms of groupoids $$\xymatrix{M \ar [r]
^{\eta} \ar [dr] _{\mu} & H \ar [d] ^{\phi}
\\ & P} $$
and an action of  $P$  on  $M$  such that  $\mu$   is a crossed
module, $\eta$    is an inclusion of a totally intransitive
subgroupoid, and if  $m \in  M , h \in  H$ and $h^{-1} mh$  is
defined, then  $h^{-1}mh = m^{\phi h}$.  It follows that  $M$   is
normal in $H$. Conversely, given such a morphism of crossed modules,
we obtain a double   $P$-module as in 1.1 with  $V = P$. Notice also
that if  $N = Ker\, \phi$, then $N$ operates trivially on $M$. }\qed
\end{example}
\begin{rem} This double category  also has a {\it thin structure} in
the sense that the set of quintuples of the form $ (1:u {\;}
^{\textstyle a} _{\T c} {\;} v)$ form a subdouble category of the
main double category. \qed
\end{rem}
\begin{rem}
It is shown in \cite{BS76} that from a double groupoid one can
recover two crossed modules, a kind of horizontal one and a vertical
one. However the groupoid conditions are needed to recover the
actions, and we know no way to do this in the category case. From a
double category one can recover two 2-categories, by restricting to
the subdouble categories where either the horizontal, or vertical,
edge categories are discrete. \qed
\end{rem}

\section{A new construction  of a homotopical  double
groupoid}\label{sec:newhomotopdg}

Let $X_*=(X;A,B;C)$ be  a space $X$ with three subspaces $A,B,C$
such that $C \subseteq A \cap B$. Let $RX_*$ be the space of maps
$f: I^2 \to X$ which map the edges $\partial ^\pm_1 I^2$ in
direction 1 into $B$, the edges $\partial ^\pm_2 I^2$ in direction 2
into $A$, and map all the vertices $\partial \partial I^2$ into $
C$. This is shown in the following diagram:

$$\def\labelstyle{\textstyle}
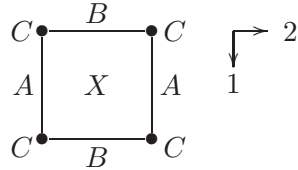
\xymatrix @=3pc @W=0pc @M=0pc
{\bullet \ar@{-}[r]^{B} \ar@{-}[d]_{A} \ar@{}[dr]|{X }
\save[]+<-0.7em,0mm>*{ C} \restore & \bullet \ar@{-}[d] ^{A}
\save[]+<0.8em,0mm>*{ C} \restore
\\ \bullet  \ar@{-}[r] _{B}\save[]+<-0.7em,-1mm>*{ C}
\restore & \save[]+<0.8em,-1mm>*{ C} \restore \bullet} \quad \xy
(0,0)*+{}="a",(0,-7)*+{\rule{0em}{1.5ex}1}="b",(7,0)*+{\;2}="c"
\ar@{->} "a";"b" \ar @{->}"a";"c" \endxy$$

The boundary maps and degeneracies give the following geometric
structure to $RX_*$, in which $RX_1$ is the set of maps $(I,
\{0,1\}) \to (A,C)$ and $RX_2$ is the set of maps $(I, \{0,1\}) \to
(B,C)$:
$$ \xymatrix{(RX_*,\circ_1,\circ_2) \ar  [d] \ar @<1ex> [d] \ar  [r] \ar @<1ex> [r]&\ar @<1ex>[l]
 (RX_2,\circ_2)\ar  [d] \ar @<1ex> [d]\\
(RX_1,\circ_1)\ar  [r] \ar @<1ex> [r]\ar @<1ex>[u]& C \ar
@<1ex>[l]\ar @<1ex>[u] }$$ Clearly the set $RX_*$ obtains two
compositions $\circ_1, \circ_2$ from the usual composition of
squares in the two directions, while $RX_2,RX_1$ have just one
composition. These compositions are of course not associative, nor
do they have identities, but they do have reverses, $-_1,-_2$. They
do however satisfy the interchange law. Further the face and
degeneracies respect the compositions. The following theorem
generalises results from \cite{BH78:sec}.
\begin{thm}
Let $\rho X_*$ be the quotient of $RX_*$ by the  relation of
homotopy rel vertices of $I^2$ and through the elements of $RX_*$.
Then the compositions $\circ_1, \circ_2$ are inherited by $\rho X_*$
to give it the structure of double groupoid if the following
condition holds:

\noindent (Con) For all $c \in C$, the induced morphisms $\pi_2(A,c)
\to \pi_2(X,c), \pi_2(B,c) \to \pi_2(X,c) $ have the same  image.

Further, under this condition, the natural morphisms $$\pi_2(X,A;C)
\to (\rho X_*,\circ_1), \quad \pi_2(X,B;C) \to (\rho X_*,\circ_2)$$
are injective.
\end{thm}

\begin{proof} The proof is a small elaboration of a similar proof for the case $A=B$ in
\cite{BH78:sec}. Details are also given in \cite{Brown-grenoble}.

The class in $\rho X_*$ of an element $\alpha$ of $RX_*$ is written
$\brclass{\alpha}$.

We develop only the horizontal case; the other follows by symmetry.
So, let us consider two elements $\langle \! \langle\alpha\rangle \!
\rangle, \langle \! \langle\beta\rangle \! \rangle \in D$ such that
$\langle \! \langle\partial^+_2\alpha\rangle \! \rangle = \langle \!
\langle\partial^-_2\beta\rangle \! \rangle$, i.e. we have continuous
maps
$$\alpha, \beta : (I^2; \partial^\pm_2 I^2, \partial^\pm_1 I^2; \partial ^2I^2) \to (X;A,B;C)$$
and a homotopy
$$ h : (I, \partial (I)) \times I \to (A,C)$$
from $\alpha|_{\{1\}\times I}$ to $\beta|_{\{0\}\times I}$ rel
vertices, i.e. $h({0} \times I) = y$ and $h({1} \times I) = x$. We
define now the composition by
$$ \langle \! \langle\alpha\rangle \! \rangle +_2 \langle \!
\langle\beta\rangle \! \rangle = \langle \! \langle \alpha +_2 h +_2
\beta \rangle \! \rangle
 = \langle \! \langle[\alpha, h, \beta]\rangle \! \rangle.$$
This is given in a diagram by

\begin{equation*}
\label{eq:2hcomp}
\def\labelstyle{\textstyle}
\vcenter{\xymatrix @=3pc @W=0pc @M=0pc { \ar@{-}[r] ^{B} \ar@{-}[d]
_{A} \ar@{}[dr]|\alpha & \ar@{-}[r] ^{c} \ar@{-}[d]|A \ar@{}[dr]|h &
\ar@{-}[r] ^{B} \ar@{-}[d]|A \ar@{}[dr]|\beta & \ar@{-}[d] ^{A}
\\ \ar@{-}[r] _{B} & \ar@{-}[r] _{d} & \ar@{-}[r] _{B} &
}}
\end{equation*}

To prove this is independent of the choices made we chose two other
representatives $\alpha' \in \langle \! \langle\alpha\rangle \!
\rangle$ and $\beta' \in \langle \! \langle\beta\rangle \! \rangle$
and a homotopy $h'$ from $\alpha'|_{\{1\}\times I}$ to
$\beta'|_{\{0\}\times I}$. Using them, we get
\begin{equation*}
\label{eq:2hcomp}
\def\labelstyle{\textstyle}
\vcenter{\xymatrix @=3pc @W=0pc @M=0pc { \ar@{-}[r] ^{B} \ar@{-}[d]
_{A} \ar@{}[dr]|{\alpha'} & \ar@{-}[r] ^{c} \ar@{-}[d]|A
\ar@{}[dr]|{h'} & \ar@{-}[r] ^{B} \ar@{-}[d]|A \ar@{}[dr]|{\beta'} &
\ar@{-}[d] ^{A}
\\ \ar@{-}[r] _{B} & \ar@{-}[r] _{d} & \ar@{-}[r] _{B} &
}}
\end{equation*}
which should give the same composition in $\rho X_*$. Let $\phi:
\alpha\simeq\alpha', \psi: \beta: \simeq \beta'$ be homotopies of
the required type. They with $h,h'$ give rise to a diagram of the
following kind.

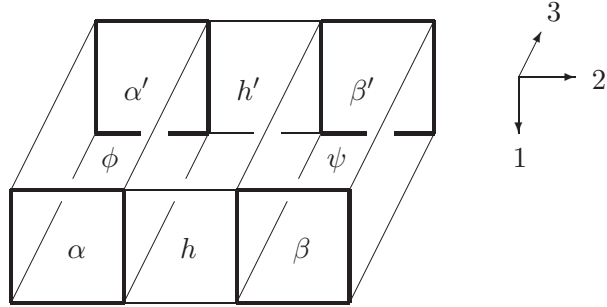
\begin{figure}[htbp]
\begin{center}
\setlength{\unitlength}{1.5cm}
\begingroup\makeatletter\ifx\SetFigFont\undefined
% extract first six characters in \fmtname
\def\x#1#2#3#4#5#6#7\relax{\def\x{#1#2#3#4#5#6}}%
\expandafter\x\fmtname xxxxxx\relax \def\y{splain}%
\ifx\x\y   % LaTeX or SliTeX?
\gdef\SetFigFont#1#2#3{%
  \ifnum #1<17\tiny\else \ifnum #1<20\small\else
  \ifnum #1<24\normalsize\else \ifnum #1<29\large\else
  \ifnum #1<34\Large\else \ifnum #1<41\LARGE\else
     \huge\fi\fi\fi\fi\fi\fi
  \csname #3\endcsname}%
\else \gdef\SetFigFont#1#2#3{\begingroup
  \count@#1\relax \ifnum 25<\count@\count@25\fi
  \def\x{\endgroup\@setsize\SetFigFont{#2pt}}%
  \expandafter\x
    \csname \romannumeral\the\count@ pt\expandafter\endcsname
    \csname @\romannumeral\the\count@ pt\endcsname
  \csname #3\endcsname}%
\fi \fi\endgroup {\renewcommand{\dashlinestretch}{30}
\begin{picture}(3.5,2.5)
\linethickness{0.4mm} \put(0,0){\line(1,0){1}}
\put(2,0){\line(1,0){1}} \put(0,0){\line(0,1){1}}
\put(0,1){\line(1,0){1}} \put(2,1){\line(1,0){1}}
\put(1,1){\line(0,-1){1}} \put(2,1){\line(0,-1){1}}
\put(3,0){\line(0,1){1}} \put(0.75,2.5){\line(1,0){1}}
\put(2.75,2.5){\line(1,0){1}} \put(0.75,2.5){\line(0,-1){1}}
\put(1.75,2.5){\line(0,-1){1}} \put(2.75,2.5){\line(0,-1){1}}
\put(3.75,2.5){\line(0,-1){1}}
% thick back lines
\put(0.75,1.5){\line(1,0){0.4}}
 \put(1.75,1.5){\line(-1,0){0.35}}
 \put(2.75,1.5){\line(1,0){0.4}}
\put(3.75,1.5){\line(-1,0){0.35}}
 \linethickness{0.1mm}
\put(1.75,2.5){\line(1,0){1}}
 \put(1,0){\line(1,0){1}} \put(1,1){\line(1,0){1}}
\put(1.75,1.5){\line(1,0){0.4}}
 \put(2.75,1.5){\line(-1,0){0.35}}
 \put(0,0){\line(1,2){0.45}} \put(0.75,1.5){\line(-1,-2){0.20}}
 \put(1,0){\line(1,2){0.45}} \put(1.75,1.5){\line(-1,-2){0.20}}
 \put(2,0){\line(1,2){0.45}} \put(2.75,1.5){\line(-1,-2){0.20}}
 \put(0,1){\line(1,2){0.75}}
\put(1,1){\line(1,2){0.75}} \put(2,1){\line(1,2){0.75}}
 \put(3,0){\line(1,2){0.75}}
 \put(3,1){\line(1,2){0.75}}
\put(0.5,0.4){\shortstack{$\alpha$}}
\put(1.5,0.4){\shortstack{$h$}}
\put(2.5,0.4){\shortstack{$\beta$}}
\put(1,1.8){\shortstack{$\alpha'$}} \put(2,1.8){\shortstack{$h'$}}
\put(3,1.8){\shortstack{$\beta'$}}
\put(0.8,1.2){\shortstack{$\phi$}}
\put(2.8,1.2){\shortstack{$\psi$}} \put(4.5,2){\vector(1,0){0.5}}
\put(4.5,2){\vector(0,-1){0.5}} \put(4.5,2){\vector(1,2){0.2}}
 \put(4.45,1.2){\shortstack{$1$}}
\put(5.15,1.90){\shortstack{$2$}}
 \put(4.75,2.5){\shortstack{$3$}}
\end{picture}
}
\caption{Filling the hole in the middle} \label{fig:Fig61}
\end{center}
\end{figure}

\noindent We seem to have a hole in the middle. The key point is
that all homotopies are rel vertices. So the bottom face of this
hole may be filled by a constant homotopy. Then we can use a
retraction to fill the hole, and this will give a cube in $A$, whose
top face is a map $(I^2, \partial \partial I^2) \to (A,c)$. By our
assumption (Con), this is deformable rel boundary and in $X$ to a
map $(I^2,
\partial
\partial I^2) \to (B,c)$. This homotopy is now added in direction 1
to the homotopy of the middle hole, and squashed down to give
another filler of the hole;  the composition in direction 2 of these
three cubes is now a homotopy rel vertices through maps $(I^2;
\partial^\pm_2 I^2, \partial^\pm_1 I^2; \partial ^2I^2) \to
(X;A,B;C)$ as required.

The verification of the groupoid axioms is entirely analogous to the
case of the fundamental groupoid. We now verify the interchange law.

Suppose given an array of composable elements of $\rho X_*$:
\begin{equation*}
  \begin{pmatrix}
\brclass{\alpha} & \brclass{\beta}\\
\brclass{\gamma}& \brclass{\delta}
  \end{pmatrix}
\end{equation*}
This gives rise to a partially filled array
\begin{equation*}
  \begin{pmatrix}
    \alpha & h & \beta \\
    k & & k' \\
    \gamma & h' & \delta
  \end{pmatrix}
\end{equation*}
However because of the rel vertices hypothesis on the homotopies
$h,h',k,k'$ the hole in the middle can be filled with a constant
map. Reading the resulting matrix in two ways gives the required
interchange law.

To this end we use the {\it connections} $\Gamma^\pm_i$ which are
available in the  cubical singular set $S^\square(X)$ of a space, as
in \cite{BH81:algcub,ABS,Grandis-cubsite}, for example. The main
point is that a connection $\Gamma:S^\square(X)_n \to
S^\square(X)_{n+1}$, defined using the functions $\max, \min$, gives
a kind of degeneracy in which two adjacent faces of $\Gamma (f)$
coincide. It is convenient to represent these symbolically as $\bl\,
, \br\, , \tl\, , \tr\,$. The traditional cubical degeneracies are
analogously represented by $\vv \, , \hh\,$. In our current
situation we surround the homotopy $H$ by connections and constant
homotopies, and also using the hypothesis (Conn) to obtain another
homotopy from this time $\alpha$ surrounded by constant maps or
$\vv$. In particular, the hypothesis (Conn) is used twice to obtain
homotopies $\xi, \xi'$ as part of the following picture of the new
homotopy.  To show this the following picture gives the picture at
$t \in [0,1]$, but the connections are actually applied on
2-dimensional faces in directions 1 and 2 of $H$. The wiggly lines
denote constant homotopies. This also illustrates that one of the
aims of bringing in connections and 2-dimensional rewriting was to
give a more algebraic method of constructing homotopies than
previously available. For another application of such rewriting,  to
rotations, see \cite{Brown-higherdimgroup}.

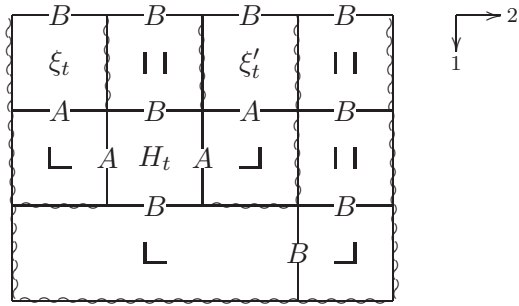
$$\xybiglabels\xymatrix@M=0pt@=3pc{\ar@{-}[ddd]\ar @{~}[ddd]\ar @{-} [r]|B
\ar@{}[dr] |{\xi_t} &\ar @{~}[d]\ar@{}[dr] |{\vv}\ar @{-}
[r]|B\ar@{-}[d]&\ar@{}[dr] |{\xi'_t}\ar @{~}[d] \ar @{-}
[r]|B\ar@{-}[d]&\ar@{-}[d]\ar @{-} [r]|B\ar@{}[dr] |{\vv}\ar
@{~}[d]&\ar@{-} [d]\ar@{~}[ddd]
\\\ar@{}[dr] |{\bl}\ar @{-}[r]|A &\ar@{}[dr] |{H_t}\ar@{-}[d]|A\ar @{-} [r]|B&
\ar@{-}[d]|A\ar@{}[dr] |{\br}\ar @{-} [r]|A&\ar@{-}[d] \ar @{-}
[r]|B\ar@{}[dr] |{\vv}\ar
@{~}[d] &\ar@{-}[d]\\
\ar @{-} [r]\ar @{~}[r]&\ar@{}[dr] |{\bl}\ar @{-} [r]|B&\ar @{-}
[r]\ar @{~}[r]&\ar@{-}[d]|B\ar @{-} [r]|B
\ar@{}[dr] |{\br}&\ar@{-}[d] \\
\ar @{-} [rrrr]\ar@{~}[rrrr] &&&&}\qquad \xdirects{2}{1}$$
\end{proof}
\begin{rem}
Even in the case $C$ is a singleton, the condition (Conn) is a non
trivial condition needed to make $\rho X_* =\pi_0(RX_*)$ a double
groupoid. Of course it is satisfied if $A=B$, giving the double
groupoid used in \cite{BH78:sec} to prove a 2-dimensional van Kampen
theorem. However it is proved in \cite{Lod82}, see also \cite{Gi},
that even without this condition the compositions $\circ_1, \circ_2$
are inherited by the group $\pi_1(RX_*)$ to give this group the
structure of cat$^2$-group, i.e. a double groupoid internal to the
category of groups.  This fact is the foundation for the work of
\cite{Lod82,BL87}. \qed
\end{rem}

\end{document}